\documentclass[11pt]{article}
\usepackage{amssymb}
\usepackage{latexsym}

\newtheorem{thm}{Theorem}[section]
\newtheorem{lem}[thm]{Lemma}
\newtheorem{cor}[thm]{Corollary}
\newtheorem{pro}[thm]{Proposition}

\newtheorem{defi}[thm]{Definition}

\newcommand{\gm }{\Gamma }
\newcommand{\lon }{\longrightarrow }
\newcommand{\be }{\begin{eqnarray*}}
\newcommand{\ee }{\end{eqnarray*}}

\newcommand{\tu}{\tilde{u}}
\newcommand{\tv}{\tilde{v}}

\newcommand{\pf}{\noindent{\bf Proof.}\ }
\newcommand{\qed}{\begin{flushright} $\Box$\ \ \ \ \ \
                  \end{flushright}}
\newcommand{\complex}{{\mathbb C}}
\newcommand{\reals}{{\mathbb R}}

\newcommand{\hstar}{*_{\hbar}}

\newcommand{\torus}{{\mathbb T}}
\newcommand{\integers}{{\mathbb Z}}

\newcommand{\half}{\frac{1}{2}}

\newcommand{\cala}{{\cal A}}

\newcommand{\calf}{{\cal F}}

\newcommand{\calj}{{\cal J}}

\newcommand{\call}{{\cal L}}

\newcommand{\calx}{{\cal X}}

\newcommand{\smalcirc}{\mbox{\tiny{$\circ $}}}

\def\description label#1{\hfil\bf[#1]\hfil}
\newcommand{\diff}[1]{\frac{d}{d#1}}
\newcommand{\pttt}[2]{\frac{\partial #1}{ \partial #2}}

\newcommand{\ih}{\frac{i}{\hbar}}
\newcommand{\hi}{\frac{\hbar}{i}}

\newcommand{\Hoch}{Hochschild }
\newcommand{\diffh}{\diff{\hbar}}
\newcommand{\td}{\tilde{D}}
\newcommand{\ff}{\calf}
\newcommand{\tphi}{\tilde{\phi }}
\newcommand{\delt}{B}
\newcommand{\lauh}{[\hbar^{-1} ,\hbar ]]}

\def\sdp{\mathbin{\hbox{$\mapstochar\kern-.3333em\times$}}}
\def\pds{\mathbin{\hbox{$\times\kern-.55em\mapstochar\,$}}}

\newcommand{\wed}{\mathbin{\lower1.5pt\hbox{$\scriptstyle{\wedge}$}}}

\let\Tilde=\widetilde

\def\chigh{{\raise1.5pt\hbox{$\chi$}}}
\let\phi=\varphi
\def\til0{\Tilde{0}}

\def\dminus{\raise2pt\hbox{\vrule height1pt width 2ex}\hskip3pt}

\def\pback#1{\mathbin{{{\lower1.2ex\hbox{$\times$}}\atop #1}}}

\def\vlra{\hbox{$\,-\!\!\!-\!\!\!-\!\!\!-\!\!\!-\!\!\!
-\!\!\!-\!\!\!-\!\!\!-\!\!\!-\!\!\!\longrightarrow\,$}}

\def\gpd{\,\lower1pt\hbox{$\longrightarrow$}\hskip-.24in\raise2pt
             \hbox{$\longrightarrow$}\,}

\def\lgpd{\,\lower1pt\hbox{$\vlra$}\hskip-1.02in\raise2pt\hbox{$\vlra$}\,}

\def\llgpd{\,\lower1pt\hbox{$\vvlra$}\hskip-1.3in\raise2pt\hbox{$\vvlra$}\,}

\begin{document}

\title{{\bf Hochschild cohomology and characteristic classes for star-products }\\Dedicated to
V.I. Arnol'd for his 60th birthday}

\author{ALAN WEINSTEIN \thanks{Research partially supported by NSF
        grant DMS-9625122, the Japan Society for the Promotion of
        Science, and the Miller Institute for Basic Research
        in Science.}\\
        Department of  Mathematics\\
        University of California\\
        Berkeley, CA 94720, USA\\
        {\sf email: alanw@math.berkeley.edu}\\
        and \\
 PING XU \thanks{ Research partially supported by NSF
        grant DMS 9704391, and Kyoto University, Japan.}\\
 Department of Mathematics\\The  Pennsylvania State University \\
University Park, PA 16802, USA\\
        {\sf email: ping@math.psu.edu }}

\maketitle

\begin{abstract}
We show that the Hochschild cohomology of the algebra obtained by formal
deformation quantization on a symplectic manifold is isomorphic to the
formal series with coefficients in the de Rham cohomology of the
manifold.  The cohomology class obtained by differentiating the
star-product with respect to the deformation parameter is seen to be
closely related to the characteristic class of the quantization.  A
fundamental role in the analysis is played by ``quantum Liouville
operators,'' which rescale the deformation parameter in the same way
in which Liouville vector fields scale the Poisson structure (or the
units of action).  Several examples are given.
\end{abstract}

\section{Introduction}

Deformation quantizations of a symplectic manifold $M$ have been
completely classified through the combined efforts of many people,
beginning with the seminal paper \cite{BFFLS} which set out the
problem and established the importance of the second de Rham
cohomology.  The role of this cohomology was further clarified in
\cite{Gu}.  The existence of deformation quantizations was
established in \cite{de-le:existence}, so that the moduli space of
deformation quantizations could then be identified as an affine space
over the vector space $H^2(M)[[\hbar]]$.   Identification of an origin
in this affine space, i.e. the choice of a distinguished isomorphism
class of quantizations, was carried out by several people, so that there
are now many versions of the classification
\cite{be-ca-gu:}\cite{Deligne}\cite{de-le:revisited}\cite{F1}\cite{NT1}.
The result is that the isomorphism type of a deformation quantization
is encoded by a {\em characteristic class} in $H^2(M)[[\hbar]]$.

In much of the work cited above, as well as the influential paper
\cite{om-ma-yo:weyl}, an important role was played by local
transformations which rescale the deformation parameter $\hbar$, as well as
the corresponding infinitesimal transformations, which we will call
{\em quantum Liouville operators}, since their classical limits are
the Liouville vector fields, whose flows rescale the symplectic
structure.  In particular, De Wilde   \cite{de:deformations}
identifies all but one term of
the characteristic class of a quantization
as the obstruction to the existence of a
quantum Liouville operator when the underlying symplectic structure is
exact (i.e. admits a Liouville vector field), while Deligne \cite{Deligne}
treats the general case and shows that the characteristic class (again
with the exception of a ``stray'' term) is the
obstruction to the existence of an object (``gerbe'') which globalizes
the local quantum Liouville operators in a somewhat more complicated way.

The purpose of the present paper is to obtain some of the results
described above by considering the derivative of a star-product with
respect to the deformation parameter as a Hochschild 2-cocycle $c$ on
the noncommutative algebra\footnote{Throughout this paper,
smooth functions and de Rham cohomology will be considered to have
{\em complex} values.}
$C^{\infty}(M)[[\hbar]]$ given by the quantization, and then attaching
to this cocycle a de Rham or \v{C}ech cocycle which represents (a
derivative of) the characteristic class.

Our first step is thus to construct a mapping
$\Phi$ (which turns out to be an isomorphism) from the Hochschild
cohomology of the deformed algebra with the \v{C}ech or de Rham
cohomology of $M$ (with a formal variable adjoined).  Although results
very close to this can be found in the literature
(e.g. \cite{br-ge:homology} and \cite{NT1}), we have not found a
complete proof of what we need.  In any case, our very simple proof,
which is a near-replica of A. Weil's proof \cite{Weil} of the
isomorphism between \v{C}ech and de Rham cohomology, seems to be new.

The second task is to compute the image $\Phi (c)$ in the cohomology
of $M$.  For this purpose, we use the method of quantum exponential
mappings \cite{Xu} to transfer the problem to the bundle of Weyl
algebras on the tangent spaces of $M$ and thereby to relate $\Phi (c)$
to the Weyl curvature of a Fedosov connection and hence to the usual
characteristic class of a deformation.  The relation matches that
already obtained in Proposition 4.4. of \cite{Deligne} by a very
different method.

In the last section, we present several examples of symplectic (and
Poisson) manifolds and their quantizations to illustrate possible
relations between the existence of classical and quantum Liouville
operators.  

On the day this manuscript was completed, we received the preprint
of Kontsevich \cite{ko:deformation}, in which the classification of
deformation quantizations is shown to be completely equivalent to the
formal classification of Poisson structures.  We hope 
that our rather direct approach to the characteristic class will
be useful in the interpretation of this equivalence.

\noindent
{\bf Acknowledgments}  We would like to thank RIMS and our host
professors Takahiro Kawai and  Kyoji Saito
for their hospitality while work on this project was being done.
Thanks go also to  Jean-Luc Brylinski, Pierre Deligne, Curt McMullen,
Marina Ratner, and Boris
Tsygan for helpful discussions.

\section{Deformation quantization}

In this section, we will recall some basic
ingredients  of
Fedosov's construction of
star-products on a symplectic manifold, as well as some
useful notation.  For details, readers should
consult \cite{F1} \cite{F2} \cite{Xu}.

Let $(M, \omega )$ be a symplectic manifold of dimension $2n$.
Each tangent space $T_{x}M$ is equipped with  a linear symplectic
structure, which can be quantized by the standard Moyal-Weyl product to
produce an algebra $W_{x}$.
More precisely,

\begin{defi}
The formal Weyl  algebra $W_{x}$  associated to $T_{x}M$
is an associative algebra  with a unit over $\complex $,
whose elements consist of formal power series in the formal parameter
$\hbar $ with coefficients being formal  polynomials
on $T_{x}M$. In other words,
each element has the form
\begin{equation}
\label{eq:general}
a(y, \hbar )=\sum \hbar^{k}a_{k, \alpha }y^{\alpha }
\end{equation}
where   $y=(y^{1}, \cdots , y^{2n})$ are linear coordinates
on $T_{x}M$, $\alpha =(\alpha_{1}, \cdots , \alpha_{2n})$
is a multi-index and         $y^{\alpha}=(y^{1})^{\alpha_{1}}\cdots
(y^{2n})^{\alpha_{2n}}$.
The product is defined by the Moyal-Weyl rule:

\begin{equation}
\label{eq:moyal}
a  *b=\sum_{k=0}^{\infty} (-\frac{i\hbar}{2})^{k}\frac{1}{k!}\pi^{i_{1}j_{1}}
\cdots \pi^{i_{k}j_{k}}
\frac{\partial^{k} a}{\partial y^{i_{1}}\cdots \partial y^{i_{k}}}
\frac{\partial^{k} b}{\partial y^{j_{1}}\cdots \partial y^{j_{k}}} .
\end{equation}
\end{defi}


Let $W=\cup_{x\in M}W_{x}$. Then $W$ is a bundle of algebras over $M$, called
the Weyl bundle over $M$.
 Its  space of smooth  sections $\gm W$  forms  an associative algebra
with unit under
fiberwise multiplication.
One may think of the sections of $W$ as the functions on a ``quantized
tangent bundle'' of  $M$.  As a vector space, $\gm W$ may be
identified with the infinite jets along the zero section of functions
on $TM$ (with $\hbar$
adjoined).  The algebra structure is a
deformation quantization of the Poisson structure given by fiberwise
Poisson bracket (using the constant symplectic structure on each
tangent space).  

 The center $Z (W)$  of $\gm W$ consists of sections
 not containing $y$'s and can be naturally
identified with $C^{\infty}(M)[[\hbar ]]$.

By assigning degrees to $y's$ and $\hbar$
with $deg y^{i}=1$ and $deg \hbar =2$, we obtain
a natural filtration
$$ C^{\infty}(M)[[\hbar ]]
\subset \gm (W_{1})\subset \cdots \gm (W_{i})\subset \gm (W_{i+1})\cdots
\subset \gm (W)    $$ with respect
to the total degree $2k+l$ of the terms in the series
 in Equation (\ref{eq:general}).

A differential form with values in $W$ is a section of the
bundle $W\otimes \wedge^{q}T^{*}M$, which  
can be expressed locally as
\begin{equation}
a(x, y, \hbar , dx)=\sum \hbar^{k}   a_{k, i_{1}\cdots i_{p}, j_{1}\cdots j_{q}}
y^{i_{1}}\cdots y^{i_{p}}dx^{j_{1}}\wedge \cdots \wedge dx^{j_{q}}.
\end{equation}
Here the coefficient $a_{k, i_{1}\cdots i_{p}, j_{1}\cdots j_{q}} $
is  symmetric with respect to
$i_{1} \cdots i_{p}$ and antisymmetric in $j_{1}\cdots j_{q}$.
For short,  we denote the space of
sections of the bundle $W\otimes \wedge^{q}T^{*}M$ by
$\gm W\otimes \Lambda^{q}$.

The fiberwise commutator in $\gm W$ extends to $\gm W \otimes \Lambda^{*}$ by
the standard procedure for differential forms with values in a bundle
of Lie algebras.
In addition, the usual   exterior derivative 
induces, in a way special to the present situation,
an operator $\delta$ on $W$-valued
differential forms:
\begin{equation}
\delta a=dx^i\wedge \frac{\partial a}{\partial y^i}, \ \ \ \forall a\in \gm W\otimes \Lambda^* .
\end{equation}
Alternatively, one can write:
\begin{equation}
\delta a=-[\ih \omega_{ij}y^{i}dx^{j} ,  a]. 
\end{equation}
Note that $\delta$ is an ``algebraic'' operator in that it does not
involve derivatives with respect to $x$.

Let $\nabla$ be a torsion-free symplectic connection on $M$ and 
$$\partial : \gm W \lon \gm W \otimes \Lambda^{1}$$
 its induced covariant
derivative.

Consider   a connection on $W$ of the form:
\begin{equation}
\label{eq:connection}
D=-\delta +\partial + \ih [\gamma ,\cdot \ ],
\end{equation}
with $\gamma \in \gm W \otimes \Lambda^{1}$.

Clearly,  $D$ is a derivation with respect to the  Moyal-Weyl product, i.e.,
\begin{equation}
D(a* b)=a * Db +Da*  b.
\end{equation}

 A simple calculation yields that
\begin{equation}
D^{2}a=-[\ih \Omega , a] , \ \forall a\in \gm W,
\end{equation}
where
\begin{equation}
\label{eq:curvature}
\Omega =\omega -R+\delta \gamma -\partial \gamma -\ih \gamma^{2}.
\end{equation}
Here  $R=\frac{1}{4}R_{ijkl}y^{i}y^{j}dx^{k}\wedge dx^{l}$ and
$R_{ijkl}=\omega_{im}R^{m}_{jkl}$ is the curvature tensor
of the symplectic connection.

A connection  of the form (\ref{eq:connection}) 
 is called {\em Abelian} if
$\Omega $ is a scalar 2-form, i.e., $\Omega \in \Omega^{2}(M)[[\hbar ]]$.
For an Abelian connection, the Bianchi identity  implies that
$d\Omega =D\Omega =0$, i.e., $\Omega \in Z^{2}(M)[[\hbar ]]$.
In this case,  $\Omega $ is called the {\em Weyl curvature}.

  \begin{thm} (\cite{F2})
\label{label:f1}
Let  $\nabla  $ be   any torsion free symplectic connection,
 and
$\Omega =\omega +\hbar \omega_{1} +\cdots \in
Z^{2}(M)[[\hbar ]]$  a  perturbation of the symplectic form
in the space  $Z^{2}(M)[[\hbar ]]$.
 There exists
a unique $\gamma \in \gm W_{3}\otimes \Lambda^{1}$
such that $D$,  given by Equation (\ref{eq:connection}),
 is an Abelian  connection which has Weyl curvature $\Omega$ and 
  satisfies
 $$\delta^{-1}\gamma =0.$$
\end{thm}

Such a  connection $D$ is often called a {\em Fedosov connection}.

Given a    Fedosov connection $D$,  the space of all parallel
sections $W_{D}$ automatically becomes  an  associative
algebra. 
Fedosov proved that $W_{D}$ can be naturally
identified with $C^{\infty}(M)[[\hbar ]]$, and therefore
 induces a star-product on $C^{\infty}(M)[[\hbar ]]$, which we will
call  a {\em Fedosov star-product}. More precisely, 
let $\sigma$ denote the projection from $W_{D}$ to its center
 $C^{\infty}(M)[[\hbar ]]$ defined as $\sigma (a)=a|_{y=0}$.

\begin{thm} (\cite{F2})
\label{thm:f2}
For any  $a(x, \hbar)\in C^{\infty}(M)[[\hbar ]]$ there is a unique
section $\tilde{a}\in W_{D}$ such that $\sigma (\tilde{a})=a$. Therefore,
$\sigma$ establishes an isomorphism between $W_{D}$ and
 $C^{\infty}(M)[[\hbar ]]$ as vector spaces. Moreover,
the equation:
\begin{equation}
a*_{\hbar}b=\sigma [ (\sigma^{-1} a)*(\sigma^{-1} b)], \ \ \ \forall a ,\ b\in
C^{\infty}(M)[[\hbar ]], 
\end{equation}
defines a star-product on $M$.
\end{thm}

If $\nabla$ is flat and $\Omega=\omega$, the Fedosov connection
is  simply given by $D=-\delta +\partial $. In this case, for any
$a\in C^{\infty}(M)$, $\tilde{a}=\sigma^{-1}(a)$
can be expressed explicitly as 
$$\tilde{a}=\sum_{k=0}^{\infty}\frac{1}{k!} (\partial_{i_{1}}\cdots
\partial_{i_{k}}a)y^{i_{1}}\cdots y^{i_{k}}, $$
which is just the Taylor expansion of $\exp _{x}^{*}a$ at
the origin. So the correspondence  $C^{\infty}(M)[[\hbar ]]
\lon W_{D}$ is indeed the pullback  by
the ($C^{\infty}$-jet  at the origin of the) usual exponential map.
Thus for a general Fedosov connection, it is sometimes useful to think
of 
 $\sigma^{-1}: C^{\infty}(M)[[\hbar ]]\lon \gm W$  as a quantum
exponential map.  

More precisely, 

\begin{defi}
A quantum exponential map is an $\hbar$-linear map
$\rho : \ \ C^{\infty}(M)[[\hbar ]] \lon  \gm W $ such that
\begin{enumerate}
\item $\rho (C^{\infty}(M)[[\hbar ]] )$ is a subalgebra of $\gm W$;
\item $\rho (a)|_{y=0}=a, \ \ \forall a\in C^{\infty}(M)[[\hbar ]]$;
\item $\rho (a)=a +\delta^{-1}da , \ \ \forall a\in C^{\infty}(M),  \ \ \mbox{mod } W_{2}$;
\item  $\rho (a)$ can be  expressed as a formal power series
in   $y$ and $\hbar$, with coefficients being derivatives of  $a$.
\end{enumerate}
\end{defi}

\begin{thm}
(\cite{Xu})
\label{thm:exponential}
Quantum exponential maps are equivalent to Fedosov connections.  Every
star-product, up to an isomorphism,
 admits a quantum exponential map and
is thus isomorphic to a Fedosov star-product.
\end{thm}
The final conclusion of this theorem can already be found in \cite{NT1}.



The  {\em characteristic class}
of a star-product algebra $\cala$ (isomorphic as an algebra over
$\complex[[\hbar]]$  to
$C^{\infty}(M)[[\hbar ]]$) is the class 
\begin{equation}
\label{eq:cha}
cl(\cala )= \frac{1}{\hbar}\Omega
\end{equation}
 in
$H^{2}(M)\lauh$, where  $\Omega$ is the Weyl curvature of a corresponding
Fedosov connection $D$.  (According to Fedosov \cite{F1}\cite{F2}, it
does not depend on the choice of connection.)

We end this section by introducing the following notation.
By $W^+$, we will denote the  extended Weyl bundle on $M$, i.e.,
the vector  bundle obtained by extending the ``ground ring''
$\complex[[\hbar]]$ to the field $\complex\lauh$ of formal Laurent
series.  It is obvious that any Fedosov connection
 $D$ naturally  extends to $W^+$, whose  space of parallel sections  is
$W_{D}\lauh$.

\section{\Hoch cohomology}

Let $ \hstar $ be a star-product on a
symplectic manifold $M$, and $\cala =C^{\infty}(M)[[\hbar ]]$
its star-product algebra. Consider $\tilde{\cala}=C^{\infty}(M) \lauh$,
the space of  formal Laurent series in $\hbar$ with coefficients
in $C^{\infty}(M)$. Then $\tilde{\cala}$ has an induced 
associative algebra structure.
Let $C^{k}(\tilde{\cala}, \tilde{\cala})$ be the space 
of $\hbar $-linear maps from $\tilde{\cala}\otimes \cdots  \otimes \tilde{\cala}$
to $\tilde{\cala}$, which are multi-differential operators
when  being restricted to $C^{\infty}(M)\otimes \cdots  \otimes C^{\infty}(M)$.
Clearly, each element of $C^{k}(\tilde{\cala}, \tilde{\cala})$ can be
identified with a multi-linear map $c=\sum_{l=-N}^{\infty} \hbar^{l} c_{l}$:
$$C^{\infty}(M)\otimes \cdots \otimes C^{\infty}(M) \lon C^{\infty}(M)\lauh,$$
where  each  $c_{l}$ is  a  $k$-multi-differential operator on $M$.

Define the  coboundary operator 
\begin{equation}
\label{eq:Hoch}
b:  C^{k}(\tilde{\cala}, \tilde{\cala})
\lon C^{k+1}(\tilde{\cala}, \tilde{\cala}),
\end{equation}
 as $\ih \tilde{b}$, where 
$\tilde{b}$ is the usual \Hoch coboundary operator:
\be
(\tilde{b}c)(u_{0}, \cdots , u_{k})&=&u_{0}\hstar c(u_{1}, \cdots , u_{k})+
\sum_{i=0}^{k-1}(-1)^{i+1}c(u_{0}, \cdots , u_{i}\hstar 
u_{i+1}, \cdots , u_{k})\\
&&+(-1)^{k+1}c(u_{0}, \cdots , u_{k-1}) \hstar u_{k} .
\ee
for $u_{0}, \cdots , u_{k} \in \tilde{\cala}$ and $c\in C^{k}(\tilde{\cala}, \tilde{\cala})$.

The \Hoch cohomology  $H^{*}(\cala , \cala )$ 
of the star-product algebra $\cala $ is defined
as  the cohomology group of this  complex.

\begin{thm}
\label{thm:hoch-weyl}
If the symplectic manifold  $M$ is 
contractible, then $H^{k}(\cala , \cala )=0$ for $k\geq 1$.
In addition,  we have  $H^{k}(W, W)=0$ for $k\geq 1$, where
$W$ is the formal Weyl algebra defined on the symplectic vector
space $\reals^{2n}$.

In fact, if  $c=\sum_{l=-N}^{\infty}  \hbar^{l} c_{l}$ is a $k$-cocycle, we can
always write $c=bH$ for some cochain $H$ of the form
$H=\sum_{l=-N}^{\infty} \hbar^{l}  H_{l}$.
\end{thm}
\pf We outline a proof here for the case  $k=2$. The general  case
follows by a longer version of the same argument.

Assume that $c=\sum_{k=0}^{\infty}\hbar^{k}c_{k}$ is a 
\Hoch 2-cocycle for the star-product. By $b_{0}$, we denote the 
\Hoch coboundary operator for the commutative algebra $C^{\infty}(M)$.
It follows from $b c=0$ that $b_{0} c_{0}=0$ ($b_{0} c_{0}$
is the $1/\hbar$ term of $b c$).
Therefore, we can write 
$$c_{0}=b_{0}T_{0}+c_{0}'$$
for some differential operator $T_{0}$ on $M$ and a bivector field
$c_{0}' \in \gm (\wedge^{2}TM )$.

(a) Assume that $T_{0}=0$.

By inspecting  the $\hbar^{0}$-term of $b c$, we have
\be
&&-\frac{i}{2}(\{u_{0}, c_{0}(u_{1}, u_{2})\}
-c_{0}(\{u_{0}, u_{1}\}, u_{2}) \\
 && +c_{0}(u_{0}, \{ u_{1}, u_{2}\}) -
\{c_{0}(u_{0}, u_{1}), u_{2}\})
+ b_{0} c_{1} =0
\ee
Taking the  alternating sum over all the permutations  of $u_{0}, u_{1}, u_{2}$, 
one obtains
$$[\pi , c_{0}]=0,$$
where $\pi$ is the Poisson tensor on $M$.
Since $H^{2}_{\pi}(M)\cong H^{2}(M)=0$, we may write $c_{0}=[X_{0}, \pi ]$ for
some vector field $X_{0}\in \calx (M)$.  (In the formal case, we use
the formal Poincar\'e lemma.)
On the other hand, it is simple to see that 
$bX_{0}=\half [X_{0}, \pi ] +O(\hbar )$.
Let $H_{0}=2X_{0}$. Then
$$c =bH_{0} + O(\hbar ). $$

(b).  In general, we have $b T_{0}=\ih b_{0}T_{0} + O(1)$.
Thus, $c-b(-i\hbar T_{0})=c_{0}-b_{0}T_{0}+O(\hbar )=c_{0}'+O(\hbar )$.
This would reduce to the situation in Case (a).

Finally, the conclusion follows by using the process above repeatedly. \qed

For a general symplectic manifold $M$, it turns out that
\begin{equation}
\label{eq:iso-k}
H^{k}(\cala , \cala )\cong H^{k}(M)\lauh .
\end{equation}
In what follows, we will prove this fact
for the case of  $k=2$ by  using Fedosov star-products. At the end, we
explain how to extend the argument to general  $k$.

\begin{thm}
\label{thm:isomorphism}
Let $W_{D}$ be  the Fedosov star-product algebra corresponding to a
Fedosov connection $D$. Then
$$H^{2}(W_{D}, W_{D})\cong H^{2}(M)\lauh.$$
\end{thm}


Our method is analogous to the proof of de Rham theorem.  The main
idea is to consider a double complex whose cohomology of horizontal
and vertical complexes are all zero except for the first row and the
first column, where the cohomologies are \Hoch cohomology and de Rham
cohomology, respectively.  Then a standard argument using spectral
sequences will lead to the conclusion.  However, below we will use a
more elementary method similar to Weil's proof \cite{Weil} of de
Rham's theorem.  This will allow us to obtain an explicit construction
of the isomorphism, which will be useful in our computation in next
section.

Consider the double complex $C^{p,q}$, where $C^{p, q}$
 is  the  tensor product  of the $p$-multilinear bundle maps
from $W^{+}$ to itself tensored with $q$-forms on $M$. 
The first  differential  $\delt : C^{p, q}\lon C^{p+1, q}$ 
is the fiberwise Hochschild coboundary operator $b$ (multiplied by the
factor $\ih$ as in Equation (\ref{eq:Hoch})), tensored with the identity 
on differential forms.  The second differential 
$\td :C^{p, q}\lon C^{p, q+1}$ is the Fedosov 
connection $D$, extended from the Weyl bundle to its ``associated bundle" 
of multilinear maps in the natural way, i.e.:
given any $\phi \in C^{p, q}$, $\td \phi \in C^{p, q+1}$
is defined by,
\be
(\td \phi)(u_{1}, \cdots , u_{p})=
D\phi (u_{1}, \cdots , u_{p})
-  \sum_{i=1}^{p} \phi (u_{1}, \cdots , Du_{i}, \cdots , u_{p}),
\ee
for any $u_{1}, \cdots , u_{p}\in W^{+}$.  It is  easy to check that
$\delt^{2}=0$ and $\td^{2}=0$. 
By $ H^{p, q}_{\delt }$ and $H^{p, q}_{\td}$, we denote the cohomology
group of the horizontal complex $\delt : C^{p, q}\lon C^{p+1, q}$,
 and  the vertical complex $\td :C^{p, q}\lon C^{p, q+1}$,  respectively.

\begin{lem}
\label{lem:double}
\begin{enumerate}
\item $\td \delt =\delt \td$;
\item $H_{\delt}^{p, q}=0$, if $p\geq 1$;
\item $H_{\td}^{p, q}=0$, if $q\geq 1$.
\end{enumerate}
\end{lem}
\pf (i) can be  proved by a straightforward verification,
 and (ii) follows from  Theorem  \ref{thm:hoch-weyl} directly.

For (iii), assume that $\phi\in C^{p, q}$ satisfying $\td \phi =0$.
It thus follows that for any $u_{1}, u_{2}, \cdots , u_{p}\in
W_{D}$, $D\phi (u_{1}, u_{2}, \cdots , u_{p})=0$.

Let $\psi  \in \gm W^{+} \otimes \Omega^{q-1}$ be the solution
of the following equation:

\begin{equation}
\label{eq:k}
\left\{ \begin{array}{ll}
D\psi &= \phi (u_{1}, u_{2}, \cdots , u_{p})\\
\psi  |_{y=0}&=0 .
\end{array}
\right.
\end{equation}

According to Theorem 5.2.5 in \cite{F2},
 this equation always has a  solution.
Moreover, $\psi$ can be made to depend on $u_{i}$ multilinearly. 
 Thus $\psi $ is a multilinear
map from $W_{D}$ to $W^+$. It is clear that $\psi$ is a local map.
It  thus extends to a $p$-multilinear bundle map from
$W$ to $W^{+}$ tensored with $(q-1)$-forms. In other words,
it is  an element in $C^{p, q-1}$,
which will be denoted by the same symbol $\psi$. Then, we have $\td \psi =\phi$.
To see this, we only need to note that both $\td \psi$ and $\phi$
are elements in $C^{p, q}$ and they coincide when being restricted
 to sections in $W_{D}$. \qed

Denote by $F^{p, q}$ the subspace of $C^{p, q}$ consisting of  those
elements $\phi$ satisfying $\td \delt  \phi =0$, and let $\ff^{p, q}\subseteq
F^{p, q}$ be  the subspace spanned by those elements which satisfy either
 the equation $\td \phi =0$ or
$\delt  \phi =0$.

\begin{pro}
\label{pro:pq}
\begin{equation}
F^{1, 0}/\ff^{1, 0} \cong  F^{0, 1}/\ff^{0, 1} .
\end{equation}
\end{pro}
\pf Given any $\phi \in F^{1, 0}$,  consider  the following  equation  for
$\psi$:
\begin{equation}
\label{eq:pq}
\delt \psi  = \td  \phi 
\end{equation}
By Lemma  \ref{lem:double} (i) and (ii), this equation always has  a
solution.
Moreover any two solutions differ by an element of $\ff^{0, 1}$.
Thus we obtain a well defined linear map 
\be
\mu : F^{1, 0} &\lon &F^{0 , 1}/\ff^{0, 1} \\
 \phi  &\lon & [\psi ]
\ee

(a).  $\mu$ is surjective.

For any $\psi \in F^{0 , 1}$, it follows from the equation
$\td \delt \psi =0$ and Lemma \ref{lem:double} (i) and (iii)
that there is $\phi \in C^{1, 0}$ such that
$\td \phi = \delt \psi$. Since $\td\delt \phi=\delt\td\phi =\delt^2 \psi=0$,  $\phi \in F^{1, 0}$.

(b).  $\ker \mu =\ff^{1, 0}$

Assume that $\mu \phi =[\psi ]=0$. By definition, $ \psi =X+Y$ with
the properties that $\td X=0$ and $\delt  Y=0$ for some $X, Y\in
F^{0 , 1}$. Thus $\delt \psi =\delt X$. Since $\td \phi =\delt \psi$,
$\td \phi =\delt X$.

On the other hand, we can always write $X=\td X_{1}$ by Lemma \ref{lem:double}. Thus
$\td \phi =\delt \td X_{1}=\td \delt X_{1}$.  Hence,
$\td (\phi -\delt X_{1})=0$.  Write $Y_{1}=\phi -\delt X_{1}$.
Then $\phi =\delt  X_{1} + Y_{1}$, which is easily seen
to be in $\ff^{1, 0}$.

Conversely, given any $\phi \in \ff^{1, 0}$, we can write
$\phi =X+Y$ such that $\td X=0$ and $\delt Y=0$.  Again by
Lemma \ref{lem:double}, we have $Y=\delt Y_{1}$ for some
$Y_{1}\in C^{0, 0}$. It thus
follows that $\td \phi =\td Y=\td \delt Y_{1}= \delt \td Y_{1}$.
Thus $\mu \phi    =[\td Y_{1} ]=0$. 

It follows from (a) and (b) that $\mu$ descends to an
isomorphism:

\begin{equation}
\label{eq:mu}
\mu : F^{1, 0}/\ff^{1, 0} \lon F^{0, 1}/\ff^{0, 1}. 
\end{equation}
 \qed

\begin{pro}
\label{pro:tau}
$$F^{1, 0}/\ff^{1, 0}\cong H^{2}(W_{D}, W_{D}) .$$
\end{pro}
\pf  Let $\phi \in F^{1, 0}$, then $\td \delt \phi =0$.
Define $\tphi$ to be   $\delt \phi$   restricted  to
$W_{D}$.  It is simple to see that
$D \tphi (u_{1}, u_{2})=0$ if  $u_{1}, u_{2}\in W_{D}$, so
 $\tphi \in C^{2}(W_{D}, W_{D})$.  Then
$b  \tphi =\delt^{2}\phi |_{W_{D}}=0$, and so we obtain
a well-defined map:

\be
\tau : F^{1, 0} &\lon  &H^{2}(W_{D}, W_{D})\\
\phi &\lon  & [\tphi ]
\ee

(a) $\tau$ is surjective.

Given any 2-cocycle  $\tphi \in C^{2}(W_{D}, W_{D})$.  It naturally
 extends to a bundle map $W\otimes W\lon W^{+}$, which will be
 denoted by the same symbol.
Then $\tphi$ is a fiberwise two-cocycle. According to Lemma \ref{lem:double},
 there is $\phi \in C^{1,0}$ such
that  $\tphi =\delt \phi $.  On the other hand, 
 for  any $u_{1}, u_{2}\in W_{D}$,
\be
(\td \tphi )(u_{1}, u_{2})&=&\td \tphi (u_{1}, u_{2}) -\tphi (Du_{1}, u_{2})-\tphi (u_{1}, Du_{2})\\
&=&\td \tphi (u_{1}, u_{2}) \\
&=&0,
\ee
since $\tphi (u_{1}, u_{2})\in W_{D}$.
This implies  that $\td \tphi$ is  identically zero
 since $W_{D}$ spans $W$ at each fiber.
Therefore, $\td  \delt \phi =\td  \tphi =0$. Hence $\phi \in F^{1,0}$
and $\tau (\phi )= [\tphi ]$.

(b) $\ker \tau  =\ff^{1, 0}$.

Assume that $\phi \in C^{1, 0}$ such that $\delt \phi =0$. 
It follows  from its definition,  that $\tphi =0$. 
On the other hand, if  $\td  \phi =0$, then $\phi (W_{D})\subset W^{+}_{D}$.
This shows that $\phi \in C^{1}(W_{D}, W_{D})$. Thus, $\tphi \in
B^{2}(W_{D}, W_{D})$ by definition. Therefore $\tau (\phi )= [\tphi ] =0$.
This shows that  $\ff^{1, 0}\subseteq \ker \tau  $.

Conversely, assume that $\tau (\phi )=0$. Then $\tphi \in C^{2}(W_{D}, W_{D})$
is a 2-coboundary. Write $\tphi =b\phi_{1}$, with $\phi_{1}\in C^{1}(W_{D}, W_{D})$.
Extend $\phi_{1}$ to a bundle map $W\lon W^+$, which will be denoted by
the same symbol. Then $\delt \phi =\delt  \phi_{1}$ when being restricted
to $W_D$. It thus follows that $\delt (\phi -\phi_{1} )=0$ identically.
Write $X=\phi -\phi_{1}$. Then $\phi =X+\phi_{1}$, where
$\delt X=0$ and $\td \phi_{1}=0$. The latter follows from
the fact that $\phi_{1} (W_D )\subseteq W_{D}^{+}$. 
 Hence $\ker \tau \subseteq \ff^{1, 0}$.

It follows from (a) and (b) that $\tau$ descends to an isomorphism:
\begin{equation}
\label{eq:tau}
\tau : F^{1, 0}/\ff^{1, 0}\cong H^{2}(W_{D}, W_{D}).
\end{equation}
\qed

\begin{pro}
\label{pro:lambda}
$$F^{0, 1}/\ff^{0, 1}\cong H^{2}(M)\lauh.  $$
\end{pro}
\pf  For any $\phi \in F^{0, 1}$, $\td \phi =D\phi$ by definition.
Hence $\delt D\phi =\delt \td \phi =0$. This means that $D\phi $ is in
the center of $\gm (W^{+})\otimes \Omega^{2}(M)$,
 so it is a scalar 2-form. Moreover, $d (D\phi )=D (D\phi )=
D^{2} \phi =0$. Thus $D\phi$ is  a closed 2-form in $Z^{2}(M)\lauh$.

Let $\lambda  : F^{0, 1}\lon H^{2}(M)\lauh$ be the map defined by
$\lambda  (\phi )=[D \phi ]$.

(a). $\lambda$ is surjective.

Given any closed 2-form $\theta \in Z^{2}(M)\lauh$. Let $\phi \in \gm (W^{+})\otimes \Omega^{1}(M)$
be any solution of the equation: $D\phi =\theta $. Thus, $\delt \td \phi =\delt D\phi =\delt \theta =0$ since
$\theta $ is a scalar 2-form. This means that 
$\phi \in F^{0, 1}$. It is clear that $\lambda  \phi =[D\phi ]=[\theta ]$.

(b). $\ker \lambda =\ff^{0, 1}$.

Assume that $\phi \in C^{0, 1}$  satisfying $\td \phi =0$. Then $D\phi =0$,
and therefore $\lambda  \phi =[D\phi ]=0$. 
On the other hand, if $\delt \phi =0$,
$\phi$ is a scalar 1-form. Then $D\phi =d\phi $. Hence  $\lambda \phi =[d \phi ]=0$. 
This shows that $\ff^{0, 1}\subseteq \ker \lambda  $.

Conversely, assume that $\lambda  \phi =0$. It follows, from definition, that
$D\phi =d\theta $ for some $\theta \in \Omega^{1}(M)\lauh$.
Write $\phi_{1}=\phi -\theta $. Then $D\phi_{1}=D\phi -d\theta =0$.
It is clear that $\delt \theta =0$. Hence $\phi \in \ff^{0, 1}$.
Therefore  $ \ker \lambda   \subseteq  \ff^{0, 1}$.

A combination of (a) and (b) implies that $\lambda$ descends
to an isomorphism:

\begin{equation}
\label{eq:lambda}
\lambda : F^{0, 1}/\ff^{0, 1}\lon  H^{2}(M)\lauh.
\end{equation}

This concludes the proof. \qed
{\bf Proof of Theorem \ref{thm:isomorphism}. }
This is a direct consequence of Propositions  \ref{pro:pq},  \ref{pro:tau}
and \ref{pro:lambda}. The desired isomorphism is established  by
\begin{equation}
\label{eq:phi}
\Phi =\lambda \smalcirc \mu \smalcirc \tau^{-1}: \  H^{2}(W_{D}, W_{D})\lon
 H^{2}(M)\lauh.
\end{equation}
\qed
\noindent{\bf Remark} The same argument  in Proposition \ref{pro:pq}
can be used to prove that
for $p\geq 1$ and $q\geq 0$,
$$F^{p, q}/\ff^{p, q}= F^{p-1, q+1}/\ff^{p-1, q+1}. $$
Similarly one can prove that
for any $k\geq 0$,
$$F^{k, 0}/\ff^{k, 0}\cong H^{k+1}(W_{D}, W_{D})$$
and
$$F^{0, k}/\ff^{0, k}\cong H^{k+1}(M)\lauh $$
using a similar argument as in  Proposition \ref{pro:tau} and \ref{pro:lambda}.
As a consequence one obtains the isomorphism:
$H^{k}(W_{D} , W_{D} )\cong H^{k}(M)\lauh$.
We are informed by Tsygan that he also  has   obtained a  proof
of this result \cite{Tsygan}.

The next result shows how the isomorphism $\Phi $ given by
 Equation (\ref{eq:phi}) depends on the  Fedosov connection.

\begin{thm}
\label{thm:eq}
Assume that $W_{D_{1}}$ and $W_{D_{2}}$ are isomorphic Fedosov
algebras coming from
 Fedosov connections $D_{1}$ and $D_{2}$, respectively.
Let $T: W_{D_{1}}\lon W_{D_{2}}$ be an isomorphism of star-products,
and $T_{*}: H^{2}(W_{D_{1}}, W_{D_{1}})\lon H^{2}(W_{D_{2}}, W_{D_{2}})$
its induced isomorphism. Then the following diagram:
\begin{equation}                         \label{tdvb}
\matrix{&&T_{*}&&\cr
        &H^{2}(W_{D_{1}}, W_{D_{1}})&\vlra&H^{2}(W_{D_{2}}, W_{D_{2}})&\cr
        &&&&\cr
     \Phi_{1}&\Bigg\downarrow&&\Bigg\downarrow&\Phi_{2}\cr
        &&&&\cr
        &H^{2}(M)\lauh&\vlra& H^{2}(M)\lauh&\cr
        &&id &&\cr}
\end{equation}
commutes.
\end{thm}
\pf  Since $T$ is local, $T$ extends to a bundle map $T: W\lon W$,
which is easily seen to be a fiberwise isomorphism. Hence,
$$T(1)=1.$$

The following lemma is essential to the proof of this theorem.

\begin{lem}
\label{lem:d12}
Under the same assumption as in Theorem \ref{thm:eq}, we have
\begin{equation}
D_{1}=T^{-1}D_{2}T.
\end{equation}
\end{lem}
\pf It is easy to check that $T^{-1}D_{2}T$ is also
a flat connection on the Weyl bundle since $T$ is
a bundle map. Consider $D=D_{1}-T^{-1}D_{2}T$. Then $Dfu=fDu$ for
any $f\in C^{\infty}(M)[[\hbar]] $ and $u\in \gm (W)$.
On the other hand, $Du=0$ for any $u\in W_{D_{1}}$.
Thus, $D$ is zero identically. This concludes the proof. \qed

Now $T$ induces an isomorphism $T_{*}: C^{p, q}\lon C^{p, q}$ 
by

$$(T_{*}\phi )(u_{1}, \cdots , u_{p})=T\phi (T^{-1}u_{1}, \cdots , T^{-1}u_{p} ), \ \ \ \ \forall u_{1}, \cdots , u_{p} \in W. $$

It follows from Lemma \ref{lem:d12} that 

$$\td_{2}\smalcirc T_{*}= T_{*}\smalcirc \td_{1}.$$

On the other hand, it is clear that $T_{*}$ commutes
with $\delt$:

$$\delt \smalcirc T_{*}= T_{*}\smalcirc \delt .$$

Therefore we have the following commuting  diagrams:

\begin{equation}                         
\matrix{&&T_{*}&&\cr
        &F_{1}^{1,0}/\ff_{1}^{1, 0}&\vlra& F_{2}^{1,0}/\ff_{2}^{1, 0}&\cr
        &&&&\cr
     \mu_{1}&\Bigg\downarrow&&\Bigg\downarrow&\mu_{2}\cr
        &&&&\cr
        &F_{1}^{0, 1}/\ff_{1}^{0, 1}&\vlra& F_{2}^{0, 1}/\ff_{2}^{0, 1}, &\cr
        && T_{*}  &&\cr}
\end{equation}

\begin{equation}                         
\matrix{&&T_{*}&&\cr
        &F_{1}^{1,0}/\ff_{1}^{1, 0}&\vlra& F_{2}^{1,0}/\ff_{2}^{1, 0}&\cr
        &&&&\cr
     \tau_{1}&\Bigg\downarrow&&\Bigg\downarrow&\tau_{2}\cr
        &&&&\cr
 &H^{2}(W_{D_{1}}, W_{D_{1}})&\vlra&H^{2}(W_{D_{2}}, W_{D_{2}}), &\cr
        && T_{*} &&\cr}
\end{equation}

and

\begin{equation}
\matrix{&&T_{*}&&\cr
        &F_{1}^{0, 1}/\ff_{1}^{0, 1}&\vlra& F_{2}^{0, 1 }/\ff_{2}^{0, 1}&\cr
        &&&&\cr
     \lambda_{1}&\Bigg\downarrow&&\Bigg\downarrow&\lambda_{2}\cr
        &&&&\cr
 &H^{2}(M)\lauh&\vlra& H^{2}(M)\lauh&\cr
        &&id &&\cr}
\end{equation}

The conclusion thus follows immediately. \qed

Given any star-product on $M$, the star-product algebra $\cala$ is isomorphic
to some Fedosov algebra $W_{D}$. So one obtains an isomorphism
$H^{k}(\cala , \cala )\cong H^{k}(M)\lauh$. Theorem \ref{thm:eq}
implies that this isomorphism is independent of the choice of the
Fedosov algebra, and therefore is canonical.  

It turns out that the isomorphism between Hochschild and de Rham
cohomology can also be constructed without the explicit choice of a
Fedosov connection.  To do this, we
consider the bundle $\calj M$ given by infinite jets of complex-valued
functions, with $\hbar$
adjoined as
usual.  This bundle carries a natural flat connection $D_\calj$, and the
choice of a star-product on $M$  endows it with a fiberwise star-product.
Its fiberwise  Hochschild coboundary operator  extends to differential 
forms  and commutes with the natural extension of $D_J$. 
 So we have a natural double complex, and we can use it as above
to get the Hochschild-de Rham isomorphism.  

In fact, according to \cite{NT1}, there is an isomorphism 
between $\calj M $ and the Weyl bundle $W$, under which
$D_\calj$ goes to a Fedosov connection $D$ in $W$. So the  double
complex arising from the jet bundle becomes essentially the same
in  the proof of Theorem \ref{thm:isomorphism}.

Finally, we note that it is possible to construct directly an
isomorphism between Hochschild and \v{C}ech cohomology by using a double
complex whose cochains are \v{C}ech cochains with values in Hochschild
cocycles for the deformed algebra.  One chooses a covering by open
subsets with all intersections contractible and uses the non-formal
version of Theorem \ref{thm:hoch-weyl}.
 
\section{Derivative of a star-product}

Let $*_{\hbar}$ be a star-product on a symplectic manifold $M$.
For any $f, g\in C^{\infty}(M)$, set
\begin{equation}
\label{eq:derivative}
c(f, g)=\diffh(f *_{\hbar}g) .
\end{equation}
Since $c$ is $1/\hbar$ times the infinitesimal deformation
corresponding to the rescaling of $\hbar$, it defines a \Hoch
2-cocycle on $\cala=C^{\infty}(M)[[\hbar ]]$. 

The main theorem is the following:

\begin{thm}
\label{thm:main}
Under the isomorphism $\Phi :H^{2}(\cala , \cala )\lon H^{2}(M)\lauh$
 as constructed in Theorem \ref{thm:isomorphism},
$$\Phi [c]=-i\hbar^{2}\diffh (cl (\cala )), $$
where $ cl (\cala )$ is  the characteristic class of
the star-product algebra (see Equation (\ref{eq:cha}) ).
\end{thm}

We shall divide our proof into several steps.  First note that $c$ induces a
family of 2-cocycles $\tilde{c}$ on the fibres of the Weyl bundle via
pull back by the quantum exponential map. A crucial step
is to write $\tilde{c} =B H$ for some   $\hbar$-linear
bundle map $H: W\lon W$.  Unfortunately,
it is difficult to find an explicit formula for this $H$.
On the other hand, there exists  another fiberwise 2-cocycle 
$c_{1}$ on the Weyl bundle obtained simply by taking  the
$\hbar$-derivative of the fiberwise star-product.
In general, $\tilde{c}$ and $c_{1}$ are quite different since
 the quantum exponential map itself involves $\hbar$. Nevertheless
we can obtain useful information about $\tilde{c}$ by relating it to $c_{1}$.
For $c_{1}$, it is simple to write it as a 2-coboundary:

\begin{pro}
Let $E=-\frac{i}{2} (\sum_{j}y_{j}\pttt{}{y_{j}})$ be $\frac{i}{2}$
times the Euler vector field  of the fiberwise  linear symplectic structure
on $TM$, where $(y_{1}, \cdots , y_{2n})$ are  the linear coordinates
in the fibers of $TM$. Then,
$$c_{1}=BE. $$
\end{pro}

Define a $\complex$-linear (NOT $\hbar$-linear) bundle map $\rho : W\lon W$:
\begin{equation}
\label{eq:rho}
\rho ( u )=\hi \dot{u}+Eu , \ \ \ \forall u \in W .
\end{equation}

Given any flat section $\tu\in W_{D}$, consider the following
equation for $v\in \gm (W)$:
\begin{equation}
\label{eq:v}
\left\{ \begin{array}{ll}
Dv&= D\rho (\tu ) \\
v&|_{y=0}=0 .
\end{array}
\right.
\end{equation}

Clearly, this equation has a unique solution since the right side of
Equation (\ref{eq:v})  is $D$-closed. Second, the solution $v$ depends on $ \tu$
$\hbar$-linearly. Hence, we obtain a $\hbar$-linear map
\be
H: W_{D}&\lon  &  \gm (W)\\
\tu  &\lon  &  v.
\ee
 
It is clear that $H$ extends to a $\hbar$-linear  bundle map from the Weyl
bundle $W$ to itself.
The following lemma indicates that $H$  coincides with $\rho$ 
on a special  class of parallel sections of $W_{D}$.
 
\begin{lem}
\label{lem:h-rho}
If $\tu\in W_{D}$ is a flat section such   that  $\tu|_{y=0}$ is
$\hbar$-independent, then
$$H(\tu )=\rho (\tu )$$
\end{lem}
\pf Since  $\rho (\tu )$ clearly satisfies the first part of
Equation (\ref{eq:v}),
it suffices to check that $\rho (\tu )$ satisfies
the initial condition.
 
\be
\rho (\tu )|_{y=0}&=& (\hi \dot{\tu}+E\tu)|_{y=0}\\
&=&\hi \dot{\tu}|_{y=0}\\
&=&\hi \diffh (\tu|_{y=0})\\
&=&0.
\ee
Here the last step follows from the assumption that
$\tu|_{y=0}$ is $\hbar$-independent. \qed
 
The following result is a direct consequence
of the lemma above together with Equation (\ref{eq:rho}).

\begin{pro}
\label{pro:tau-char}
 For any $f, g \in C^{\infty}(M)$,
$$c(f, g)= (BH) (\tu ,\tv)|_{y=0},$$
where  $c$ is defined as in Equation (\ref{eq:derivative}),
$\tu$ and $\tv$ are parallel sections such that $\tu|_{y=0}=f$
and $\tv|_{y=0}=g$, respectively.

In other words, $$\tilde{c}=\delt H .$$
\end{pro}
 
The next step is to write $\td H=BK_{0}$ 
for some $K_{0}\in \gm(W) \otimes \Lambda^{1}$.
 For this purpose we need the following

\begin{pro}
\label{pro:rhotu}
For any $\tu \in W_{D}$,
we have
\begin{equation}
\label{eq:rhotu}
D\rho ( \tu )=[\ih K_{0}, \tu ],
\end{equation}
where
\begin{equation}
K_{0}=-(Er-i\hbar \dot{r}+ir +\frac{i}{2}
 \omega_{ij}y^{i}dx^{j} ) \in \gm(W) \otimes \Lambda^{1}.
\end{equation}
\end{pro}

An immediate consequence is

\begin{cor}
\label{cor:rhotu}
We have 
$$\td H =BK_{0}. $$
\end{cor}
\pf Equation  (\ref{eq:rhotu})  implies that $(\td H )(\tu )=
(BK_{0})(\tu )$, if $\tu \in W_{D}$ and $\tu|_{y=0}$ is $\hbar$-independent
according to    Lemma \ref{lem:h-rho}. The conclusion thus follows
immediately since $\td H$ and $BK_{0}$ are $\hbar $-linear and 
$W_{D}$ spans each fiber of $W$.
\qed

To prove Proposition  \ref{pro:rhotu}, we need a couple of lemmas first.

\begin{lem}
\begin{enumerate}
\item $[\partial , E]=0$;
\item $[\delta , E]=-\frac{i}{2}\delta =-
[\frac{1}{2\hbar}\omega_{ij}y^{i}dx^{j} , \cdot ]$
\end{enumerate}
\end{lem}
\pf (1) can be easily verified, and is left to the
reader.

As for (2), let $a\in \gm (W)\otimes \Lambda^{|a|}$,
\be
[\delta , E]a&=&\delta Ea-E\delta a \\
&=&\delta (-\frac{i}{2}\sum_{j}y^{j} \frac{\partial a}{\partial y^{j}} )
-E(\sum_{k}dx^{k}\wedge \frac{\partial a}{\partial y^{k}})\\
&=&-\frac{i}{2}\sum_{j} dx^{j}\wedge  \frac{\partial a}{\partial y^{j}}\\
&=&-\frac{i}{2}\delta a
\ee
\qed

\begin{lem}
\label{lem:commute}
For any $a\in \gm (W)\otimes \Lambda^{|a|}$,
\begin{enumerate}
\item $[D, \diffh]a=-[(\ih r)^{\cdot}, a]-\ih (c_{1}(r, a )-(-1)^{|a|}c_{1}(a, r))$;
\item $[D, E]a= \frac{i}{2}\delta a -[\ih Er, a]+(c_{1}(r, a )-(-1)^{|a|}c_{1}(a, r))$.
\end{enumerate}
\end{lem}
\pf 
\be
&&[D, \diffh]a \\
&=&D\dot{a}-\diffh (-\delta a +\partial a +[\ih r,  a ] )\\
&=&D\dot{a}-(-\delta \dot{a}+\partial \dot{a}+[(\ih r)^{\cdot}, a]\\
&&+ [\ih r,  \dot{a}]+c_{1}(\ih r , a )-(-1)^{|a|}c_{1}(a, \ih r))\\
&=&D\dot{a}-(D\dot{a}+[(\ih r)^{\cdot}, a]+\ih c_{1}(r, a)-(-1)^{|a|} \ih c_{1}(a, r))\\
&=&-[(\ih r)^{\cdot}, a]-\ih (c_{1}(r, a )-(-1)^{|a|}c_{1}(a, r)) .
\ee

(2) $$EDa=-E\delta a +E\partial a +E[\ih r,  a ].$$

\be
E[\ih r,  a ]&=&\ih (E(r*a)-(-1)^{|a|}E(a*r))\\
&=&\ih Er *a+r* \ih Ea-c_{1}(r, a)\\
&&+(-1)^{|a|}(-\ih Ea *r-a*\ih Er +c_{1}(a, r))\\
&=&[\ih Er, a ]+[r, \ih Ea]-c_{1}(r, a)+(-1)^{|a|} c_{1}(a, r)
\ee

On the other hand, 
$$DEa=-\delta Ea +\partial E a +[\ih r,  Ea ]. $$
Hence,
\be
[D, E]a&=&-[\delta ,E]a+[\partial , E]a-[\ih Er, a]+c_{1}(r, a)-(-1)^{|a|}c_{1}(a, r)\\
&=&-[\delta ,E]a -[\ih Er, a]+c_{1}(r, a)-(-1)^{|a|}c_{1}(a, r)\\
&=& \frac{i}{2}\delta a -[\ih Er, a]+(c_{1}(r, a )-(-1)^{|a|}c_{1}(a, r)) .
\ee
\qed
{\bf Proof of Proposition \ref{pro:rhotu}}
 It follows from Lemma \ref{lem:commute} that
 $$D\dot{\tu}=-[(\ih r)^{\cdot}, \tu]-\ih (c_{1}(r, \tu  )-c_{1}(\tu , r)), $$
and
$$DE\tu= \frac{i}{2}\delta \tu -[\ih Er, \tu]+(c_{1}(r, \tu )-c_{1}(\tu , r)). $$
Thus,
\be
D\rho ( \tu )&=&\hi D\dot{\tu}+DE\tu\\
&=&-\hi[(\ih r)^{\cdot}, \tu]-[\ih Er, \tu]+\frac{i}{2}\delta \tu \\
&=&[\ih K_{0}, \tu ],
\ee
where $K_{0}=-(Er-i\hbar \dot{r}+ir +\frac{i}{2} \omega_{ij}y^{i}dx^{j} )$. \qed

Finally we need to compute $DK_{0}$:

\begin{pro}
\label{pro:dk}
\begin{equation}
DK_{0}=i(\Omega -\hbar \dot{\Omega})=-i\hbar^{2} \diffh(\frac{1}{\hbar}\Omega),
\end{equation}
where $\Omega$ is the Weyl curvature of the Fedosov connection $D$. 
\end{pro}
\pf Let $L: \gm (W)\otimes \Lambda^{k}\lon \gm (W)\otimes \Lambda^{k}$ be
the operator 
$L=E-i \hbar \diffh $.
According to Lemma \ref{lem:commute}, we have,

$$[D, \diffh]r=-[(\ih r)^{\cdot}, r]- 2\ih c_{1}(r, r ), $$
and
$$[D, E]r= \frac{i}{2}\delta r -[\ih Er, r]+2c_{1}(r, r ) .$$
Thus,
\be
[D, L]r&=&[D, E]r-i \hbar  [D, \diffh ]r\\
&=&\frac{i}{2}\delta r  -[\ih Er, r]-[\dot{r}, r]+\frac{2}{\hbar}r^{2}.
\ee
On the other hand, we have
\begin{equation}
\Omega=\omega -R -Dr+\ih r^{2} .
\end{equation}
Thus,
$$Dr=\omega- \Omega-R+\ih r^{2}.$$
Hence
$$EDr=-ER+\ih E(r\cdot r)=iR+\ih E(r\cdot r),$$
and
\be
-i \hbar \diffh Dr&=&-i \hbar (-\dot{\Omega}+(\ih r^{2})^{\cdot})\\
&=&i \hbar \dot{\Omega} -\frac{1}{\hbar}r^{2}+[\dot{r}, r]+c_{1}(r, r ).
\ee
Hence,
\be
LDr&=&EDr-i \hbar (Dr)^{\cdot}\\
&=&iR+\ih E(r\cdot r )+i  \hbar \dot{\Omega} -\frac{1}{\hbar}r^{2}
[\dot{r}, r]+c_{1}(r, r )
\ee
and
\be
DLr&=&[D, L]r+LDr\\
&=&\frac{i}{2}\delta r+\frac{1}{\hbar}r^{2}+iR +i\hbar \dot{\Omega} .
\ee

On the other hand,
\be
&&D( \frac{i}{2} \omega_{ij}y^{i}dx^{j} )\\
&=&-\frac{i}{2} \omega_{ij}dx^{i}\wedge dx^{j}+[\ih r , \frac{i}{2} \omega_{ij}y^{i}dx^{j}]\\
&=& -i \omega+\frac{i}{2}[\ih r, \omega_{ij}y^{i}dx^{j}]\\
&=&-i \omega +\frac{i}{2}[\ih \omega_{ij}y^{i}dx^{j} , r]\\
&=&-i \omega -\frac{i}{2}\delta r .
\ee
Therefore
\be
DK_{0}&=&-(DLr+iDr +D( \frac{i}{2} \omega_{ij}y^{i}dx^{j} ))\\
&=&i(\Omega -\hbar \dot{\Omega })\\
&=&-i\hbar^{2} \diffh(\frac{1}{\hbar}\Omega) .
\ee
This concludes the proof.  The cancellation of the expressions
involving $c_{1}$ still seems rather mysterious to us.
\qed 
{\bf Proof of Theorem \ref{thm:main}}
Proposition \ref{pro:tau-char} implies that $\tau^{-1}[\tilde{c}]=[H]$.
On the other hand, it follows from Corollary  \ref{cor:rhotu} 
that  $\mu [H]=[K_{0}]$. Now $\lambda [K_{0}]=[DK_{0}]=
-i\hbar^{2} [\diffh(\frac{1}{\hbar}\Omega)]$ 
according to Proposition \ref{pro:dk}.
Thus the conclusion follows immediately.  \qed

We define a {\em quantum Liouville operator} to be an $\hbar$-linear
local operator $X$ on $\cala$ such that $\hbar\diffh+X$ is a
derivation (intuitively, it generates a 1-parameter group of
algebra automorphisms which rescale $\hbar$).  It is easy to see that $X$ is
such an operator if and only if its Hochschild coboundary is a
constant multiple of the derivative cocycle $c$.  The next result then follows
immediately from Theorem \ref{thm:main}.

\begin{cor}
\label{cor:liouville}
A star-product algebra $\cala$ admits a quantum Liouville operator if and only if 
$\diffh (cl (\cala )) = 0$.
\end{cor}

\noindent{\bf Remark} This consideration of quantum
Liouville operators connects our
Theorem \ref{thm:main} to  Theorem  4.4
of  Deligne  in \cite{Deligne}.

\section{Some examples}
We have just seen that the characteristic class of a deformation
quantization is nearly the same as the obstruction to the existence of
a global Liouville operator.  Thus, it is interesting to have an
example where the characteristic class is nonzero, but a Liouville
operator exists nevertheless.
We will exhibit such an example on the cotangent bundle of a
2-torus.  

On $M=T^* \torus ^2$, with coordinates
$(\theta_1,\theta_2,p_1,p_2)$, we consider the $\hbar$-dependent
symplectic structure $\omega=\omega_0 + \hbar \omega_1$, where $\omega_0$
is the canonical symplectic structure and $\omega_1$ is the
cohomologically nontrivial form $d\theta_1\wedge d\theta_2$.  The
corresponding Poisson tensor $\pi$ is the sum of the canonical Poisson
structure $\pi_0$ and $\hbar$ times
the bivector $\pi_1 = -\frac{\partial}{\partial
p_1}\wedge \frac{\partial}{\partial p_2}$. 
 Denote by $X$ the usual
Liouville vector field $p_1\frac{\partial}{\partial p_1} + 
p_2 \frac{\partial}{\partial p_2}$. 

We obtain a star-product on $M$ by applying the usual exponential formula
for the Moyal-Weyl product (see Equation (\ref{eq:moyal}))
to the $\hbar$-dependent bivector field
$\pi$.    This is a Fedosov star-product for which $\Omega = \omega$, so
its characteristic class is $\frac{1}{\hbar}[\omega_0] +[\omega_1]$.  
Since $\omega_0$ is exact but $\omega_1$ is not, the characteristic
class is nonzero, but its derivative with respect to $\hbar$ vanishes.

It follows from Corollary \ref{cor:liouville} that this
star-product should admit a quantum Liouville operator.  In fact, it is
easy to check (using the usual rule for the derivative of an
exponential and the Lie derivative formulas $\call_X \pi_0 = \pi_0$
and $\call_X \pi_1 = 2 \pi_1$) that the Liouville vector field $X$
functions as a quantum Liouville operator for our star-product (just as
it does for the usual Moyal-Weyl product corresponding to $\omega_0$).

We note further that, replacing the term $\hbar \omega_1$ by $\hbar^2
\omega_1$ in the definition of our $\hbar$-dependent symplectic
structure, we obtain a deformation quantization for the standard
symplectic structure which does {\em not} admit a quantum Liouville operator.

Since the symplectic structure on a compact manifold cannot be
exact, no star-product on a compact symplectic manifold can admit a
quantum Liouville operator.  
The following example shows that we can 
find such an operator on a compact {\em Poisson}
manifold.  

There do exist compact regular Poisson manifolds admitting Liouville
vector fields. For instance (see \cite{we:modular} for more details),
we can take $M$ to be the unit tangent bundle of a compact surface
$\Sigma$ of
negative curvature, with a Poisson structure for which the symplectic
leaves are the unstable manifolds of the geodesic flow.  One sees from
the structure equations of the group $SL(2,\reals)$, of which $M$ is a
quotient, that the generator $X$ of the geodesic flow is a Liouville
vector field.  For such a tangential Poisson structure, the
obstruction to lifting $X$ to a quantum Liouville operator lies in the
2nd de Rham cohomology along the symplectic leaves.  At our
instigation, Ratner \cite{ra:private} has provided a proof that this
cohomology vanishes as long as $1/4$ is not an eigenvalue of the
laplacian on functions on $\Sigma$, so that the quantum Liouville field must
exist.\footnote{Ratner's proof uses the representation theory of
$SL(2,\reals)$ on the functions on $M$.  A related result for
measurable cochains which are smooth only along the leaves, based on
the nonexistence of a transverse invariant measure, is given as
Theorem 4.27 in \cite{mo-sc:global}.  Incidentally, vanishing of the
(smooth) cohomology implies that any two Poisson structures on $M$
associated with this foliation are isomorphic (at least up to
sign), as are any two tangential star-products quantizing such a
Poisson structure.}

There is an alternative way to obtain the quantum Liouville operator
for $M$; it works even if $\Sigma$ has 1/4 as an eigenvalue.
The Poisson structure on $M$
is generated by an action of a 2-dimensional solvable group $G$
of triangular matrices contained in  $SL(2,\reals)$, so it can be
quantized with the aid of a quantum  $R$-matrix, as in  \cite{Xu:1993}.
This $R$ matrix gives a left invariant quantization of a left
invariant Poisson structure on $G$, which can be given a
left-invariant quantization by Fedosov's method.  An equivariant
version of our method in this paper shows that the obstruction to
finding a left-invariant quantum Liouville operator lies in the
cohomology of left-invariant forms on $G$, which is zero.  The
resulting operator $Y$ can then be ``pushed forward'' by the
$G$ action to give a Liouville operator on $M$.

Here is an example which falls just short of having a quantum
 Liouville operator.
 Let $A$ be an element of $SL(3,\integers)$ whose
eigenvalues $\lambda_j$ are all real, with $0 < \lambda_1 < 1 < \lambda_2
\leq \lambda_3$.  Let $\pi$ be a translation invariant Poisson
structure on $\torus^3$ whose symplectic leaves are planes parallel to
the 2-dimensional expanding subspace of $A$.  The linear
transformation defined by $A$ acts on the torus, where it multiplies
the Poisson structure by $\lambda = \lambda_2 \lambda_3$; it has a
similar effect on the corresponding Moyal-Weyl product, rescaling the
deformation parameter by $\lambda$.  Thus we have a discrete family of
$\hbar$ rescaling transformations, but not a 1-parameter group.  In
fact, it follows from consideration of its fundamental homology class
(see \cite{we:modular}) that this Poisson
structure does not admit a Liouville vector field (i.e. a vector field
satisfying $\call_X \pi = \pi$).

By multiplying $\torus^3$ by $\reals$ and using a suitable symplectic
structure, one can obtain a symplectic manifold admitting no Liouville
vector field, but admitting a discrete
group of maps which scale the symplectic structure nontrivially.

Another interesting phenomenon occurs in the case of $SU(2)$ as a
Poisson Lie group with the Bruhat-Poisson structure
\cite{lu-we:poisson}. This Poisson manifold admits no Liouville vector
field (as can be seen from consideration of its modular vector field),
but it does admit a 1-parameter group of {\em continuous}
mappings which scale the Poisson structure.  (This makes sense because
the mappings are smooth on each symplectic leaf.)  This reflects the
fact \cite{na:} that the $C^*$-completions of the algebras of functions on the
corresponding quantum $SU(2)$'s are isomorphic as the deformation
parameter varies through nonzero values.

Finally, we remark that the independence of the isomorphism class of a
deformation quantization on the (nonzero) value of the deformation
parameter is important for the link between deformation quantization and
the ``E-theory'' of asymptotic morphisms of $C^*$ algebras
\cite{co-hi:deformations}.

\end{document}